\documentclass[11pt]{article}
       \usepackage{amsfonts}
       \usepackage{stmaryrd}
       \usepackage{latexsym,amssymb,mathrsfs,fancyhdr}
       \font\tenmsb=msbm10
       \font\sevenmsb=msbm7
       \font\fivemsb=msbm5
       \catcode`\@=11
       \ifx\amstexloaded@\relax\catcode`\@=\active
       \endinput\else\let\amstexloaded@\relax\fi
       \def\spaces@{\space\space\space\space\space}
       \def\spaces@@{\spaces@\spaces@\spaces@\spaces@\spaces@}
       \def\space@.  {\futurelet\space@\relax}
       \space@.   %
       \def\Err@#1{\errhelp\defaulthelp@\errmessage{AmS-TeX error: #1}}
       \def\relaxnext@{\let\next\relax}
       \def\accentfam@{7}
       \def\noaccents@{\def\accentfam@{0}}
       \def\Cal{\relaxnext@\ifmmode\let\next\Cal@\else
       \def\next{\Err@{Use \string\Cal\space only in math mode}}\fi\next}
       \def\Cal@#1{{\Cal@@{#1}}}
       \def\Cal@@#1{\noaccents@\fam\tw@#1}
       \def\Bbb{\relaxnext@\ifmmode\let\next\Bbb@\else
       \def\next{\Err@{Use \string\Bbb\space only in math mode}}\fi\next}
       \def\Bbb@#1{{\Bbb@@{#1}}}
       \def\Bbb@@#1{\noaccents@\fam\msbfam#1}
       \newfam\msbfam
       \textfont\msbfam=\tenmsb
       \scriptfont\msbfam=\sevenmsb
       \scriptscriptfont\msbfam=\fivemsb
\usepackage{dsfont}
\usepackage[german,english]{babel}
\usepackage{amsmath,amssymb}
\usepackage[square, comma, sort&compress, numbers]{natbib}
\usepackage{cases}
\usepackage{mathrsfs}
\usepackage{amsmath}
\usepackage{amsthm}
\usepackage{amsfonts}
\usepackage{amssymb}
\usepackage{latexsym}
\usepackage{fancyhdr}
\usepackage{extarrows}
\usepackage{geometry}
\newtheorem{thm}{Theorem}[section]
\newtheorem{prop}[thm]{Proposition}
\newtheorem{lem}[thm]{Lemma}
\newtheorem{rem}[thm]{Remark}
\newtheorem{iteration lemma}[thm]{iteration Lemma}
\newtheorem{cor}[thm]{Corollary}

\newtheorem{eg}[thm]{Example}

\newtheorem*{acknowledgements*}{ACKNOWLEDGEMENTS}

\begin{document}

\setlength{\columnsep}{5pt}
\title{\bf   On the g$\pi$-Hirano invertibility  in Banach algebras}
\author{\ Honglin  Zou\footnote{E-mail: honglinzou@163.com},
\ Tingting Li\footnote{E-mail: littnanjing@163.com},
\ Yujie Wei\footnote{E-mail: yujiewei1994@163.com},
\\
$^{\ast}$College of Basic Science, Zhejiang Shuren University, \\ Hangzhou 310015, China\\
$^{\dagger}$School of Mathematical Sciences, Yangzhou University, \\ Yangzhou 225002, China\\
$^{\ddagger}$School of Mathematics and Statistics, Hubei Normal University,\\ Huangshi 435002, China}
\date{}

\maketitle
\begin{quote}
{\textbf{Abstract:} \small
In a Banach algebra,  we introduce a new type of generalized inverse called g$\pi$-Hirano inverse.
Firstly, several existence criteria and the equivalent definition of this inverse are investigated.
Then, we discuss the relationship between the g$\pi$-Hirano invertibility of $a$, $b$ and that of the sum $a+b$ under some weaker conditions.
Finally, as applications to the previous additive results,
some equivalent characterizations for the g$\pi$-Hirano invertibility of the anti-triangular matrix over Banach algebras are obtained.
In particular, some results in this paper are different from the corresponding ones of classical generalized inverses, such as Drazin inverse and generalized Drazin inverse.


}
\end{quote}

\section{Introduction }\label{a}
Let $\mathcal{A}$ be a complex Banach algebra with unity 1.
For $a \in \mathcal{A}$, denote  the spectrum and the spectral radius of $a$ by $\sigma(a)$ and $r(a)$, respectively.
$\mathcal{A}^{nil}$  and $\mathcal{A}^{qnil}$ stand for
the sets of all nilpotent and quasinilpotent elements ($\sigma(a)$ =$\{0\}$) in $\mathcal{A}$, respectively.
It is well known that  $a\in \mathcal{A}^{qnil}$ if and only if $r(a)=0$.
The double commutant of an element  $a \in \mathcal{A}$ is defined by comm$^{2}$($a$)=$\{b\in \mathcal{A}: bc=cb,\ \mbox{for\ any}\ c \in \mathcal{A}\ \mbox{satisfying}\ ca=ac\}$.

As is known to all, Drazin inverse \cite {Drazin58} is a kind of classic generalized inverse and has many applications.
Until now, there have been many types of generalized inverses related to Drazin inverse.
Here we list some of them as follows.

The generalized Drazin inverse (or g-Drazin inverse)  of $a\in \mathcal{A}$ \cite{Koliha96} is the element $x\in \mathcal{A}$  which satisfies
$$ xax=x, \ \ ax=xa \ \ \mbox{and} \ \ a-a^{2}x \in \mathcal{A}^{qnil}.$$
Such $x$, if it exists, is unique and will be denoted by  $a^{d}$.

An element   $x\in \mathcal{A}$ is called the generalized strong Drazin inverse (or gs-Drazin inverse)  of $a\in \mathcal{A}$  \cite{Dijana16}
if it satisfies
$$ xax=x,\ \ ax=xa\ \   \mathrm{and} \ \ a-ax\in \mathcal{A}^{qnil}.$$

Recently, the notion of $\pi$-Hirano inverse \cite{GHS} was introduced in Banach algebras. Namely,
the $\pi$-Hirano inverse of $a\in \mathcal{A}$ is the unique element $x$ satisfying
$$ xax=x,\ \ ax=xa\ \   \mathrm{and} \ \ a-a^{n+2}x\in \mathcal{A}^{nil},$$
for some $n\in \mathbb{N}$.
Motivated by this notion, we give the  definition of the generalized $\pi$-Hirano inverse as follows.

An element   $x\in \mathcal{A}$ is called the generalized $\pi$-Hirano inverse (or g$\pi$-Hirano inverse) of $a\in \mathcal{A}$
if it satisfies
$$ xax=x,\ \ ax=xa\ \   \mathrm{and} \ \ a-a^{n+2}x\in \mathcal{A}^{qnil},$$
for some $n\in \mathbb{N}$.

All the time different types of generalized inverses were investigated in several directions
(existences, sums, block matrices, reverse order laws, applications  etc.) and in different settings (operator algebras, $C^{\ast}$-algebras, Banach algebras, rings etc.).
For example,  Drazin \cite{Drazin58} proved that $a\in R$ is Drazin invertible if and only if it is strongly $\pi$-regular (i.e. $a^{m}\in a^{m+1}R\cap Ra^{m+1}$, for some $m\in \mathbb{N}$) in a ring $R$.
Meanwhile, the Drazin invertibility of the sum $a+b$ was studied under the condition $ab=ba=0$.
Later, in a Banach algebra Koliha \cite{Koliha96} claimed  that $a\in \mathcal{A}$ has the generalized Drazin inverse if and only if 0 is not an accumulation point of  $\sigma(a)$.
For the ring case, Koliha and Patri\'{c}io \cite{KP02} showed that $a\in R$ is generalized Drazin invertible if and only if $a$ is quasipolar.
In \cite{GHS}, the authors considered the $\pi$-Hirano invertibility of a 2$\times$2 operator matrix.
More results on the generalized inverses related to this paper can be found in
\cite{CS2019F,CS2019,CS20193,C2009,Cvetkovic09,Dijana19,ZCZW2020,ZDZC19,C2017}.

All the results mentioned above served as motivation for further consideration of the g$\pi$-Hirano inverse in Banach algebras.
This paper is composed of  four sections.
In Section~2, we characterize the g$\pi$-Hirano inverse by means of the quasinilpotent elements.
Then, the equivalent definition of this inverse is given.
In Section 3, sufficient and necessary conditions for the g$\pi$-Hirano invertibility of the sum $a$+$b$ are obtained
under some weaker conditions.
In Section 4, we investigated the g$\pi$-Hirano invertibility of several kinds of  anti-triangular matrices over Banach algebras.

\medskip
Next, we introduce some  well-known lemmas, which are related to the quasinilpotency in a Banach algebra.

\begin{lem}\label{lem1} \emph{\cite[Lemma 2.1]{C2009}}
Let $a,b\in \mathcal{A}$ be such that $ab=ba$. The following hold:

$(1)$ If $a\in \mathcal{A}^{qnil}$ \emph{(}or $b\in \mathcal{A}^{qnil}$\emph{)}, then $ab\in \mathcal{A}^{qnil}$.

$(2)$ If $a,b\in \mathcal{A}^{qnil}$, then $a+b\in \mathcal{A}^{qnil}$.
\end{lem}

\begin{lem}\label{lem121} \emph{\cite[Lemma 2.4]{CK04}}
Let $a,b\in \mathcal{A}^{qnil}$.
If $ab=0$, then $a+b\in \mathcal{A}^{qnil}$.
\end{lem}

\begin{lem} \label{lem18}  \emph{\cite[Lemma 1.1]{Dijana19}}
Let $n\in \mathbb{N}$. Then,
$a\in \mathcal{A}^{qnil}$ if and only if $a^{n}\in \mathcal{A}^{qnil}$.
\end{lem}

\begin{lem}\label{lem3} \emph{\cite[Lemma 2.2]{CS2019F}}
Let $a\in \mathcal{A}$. Then,  $a$ is gs-Drazin invertible if and only if $a-a^{2}\in \mathcal{A}^{qnil}$.
\end{lem}

\section{Characterizations for the g$\pi$-Hirano invertibility }\label{a}
In this section, we investigate the existence criterion for the g$\pi$-Hirano inverse in terms of quasinilpotent elements in  Banach algebras.
Then, using this  characterization we obtain the equivalent definition for the g$\pi$-Hirano inverse.

\medskip
Firstly, we give the relationship between the g$\pi$-Hirano inverse and the g-Drazin inverse.
Let $\mathcal{A}^{d}$ and $\mathcal{A}^{g\pi H}$ denote the sets of all g-Drazin and g$\pi$-Hirano invertible elements in  $\mathcal{A}$, respectively.

\begin{prop}\label{prop1}
Let  $a\in \mathcal{A}$. If $x$ is the g$\pi$-Hirano inverse of $a$,
then $a\in \mathcal{A}^{d}$ and $a^{d}=x$.
\end{prop}
\proof Suppose that $x$ is the g$\pi$-Hirano inverse of $a\in \mathcal{A}$, i.e.
$$xax=x, ax=xa \ \mbox{and}\ a-a^{n+2}x\in \mathcal{A}^{qnil},$$
for some $n\in \mathbb{N}$.
Hence, by Lemma \ref{lem1}(1) we get
$$a-a^{2}x=(a-a^{n+2}x)(1-ax)\in \mathcal{A}^{qnil}.$$
So, $a\in \mathcal{A}^{d}$ and $a^{d}=x$.
\qed

\medskip
From Proposition \ref{prop1}, we see that the g$\pi$-Hirano inverse is a subclass of the g-Drazin inverse.
According to the uniqueness of the g-Drazin inverse, we obtain that
the g$\pi$-Hirano inverse is unique if it exists.
So, we use $a^{g\pi H}$ to denote the g$\pi$-Hirano inverse of $a$ in a Banach algebra.

\medskip
In \cite[Theorem 2.1]{{GHS}}, the authors investigated the existence of $\pi$-Hirano inverse by means of the nilpotent element.
Inspired by this theorem,  we obtain the corresponding result for the g$\pi$-Hirano inverse with the help of Lemma \ref{lem3}.

\begin{thm}\label{theorem4}
Let $a\in \mathcal{A}$.  Then the following are equivalent:

$(1)$ $a\in \mathcal{A}^{g\pi H}$;

$(2)$ $a-a^{n+1}\in \mathcal{A}^{qnil}$, for some $n\in \mathbb{N}$;

$(3)$ $a^{m}-a^{n}\in \mathcal{A}^{qnil}$, for some $m, n \in \mathbb{N}$ such that $m\neq n$.
\end{thm}
\proof  $(1) \Rightarrow (2)$.  Suppose that $a\in \mathcal{A}^{g\pi H}$.
Then, there exists $x\in \mathcal{A}$ such that
$$xax=x, ax=xa\ \mbox{and}\ a-a^{n+2}x\in \mathcal{A}^{qnil},$$
for some $n\in \mathbb{N}$.
Therefore, we have
$$a-a^{n+1}=(a-a^{n+2}x)(1+a^{n+1}x-a^{n})\in \mathcal{A}^{qnil}.$$

$(2) \Rightarrow (3)$. It is obvious.

$(3) \Rightarrow (1)$. Suppose that $n>m$.
Note that
$$\begin{array}{ccl}
$$& &(a-a^{n-m+1})^{m}=\left(a(1-a^{n-m})\right)^{m}=a^{m}(1-a^{n-m})(1-a^{n-m})^{m-1}\\
& & \phantom{(a-a^{n-m+1})^{m}}=(a^{m}-a^{n})(1-a^{n-m})^{m-1}
\end{array}$$
and $a^{m}-a^{n}\in \mathcal{A}^{qnil}$.
So, by Lemma \ref{lem1}(1) and Lemma \ref{lem18} we get $a-a^{n-m+1}\in \mathcal{A}^{qnil}$,
which gives $a^{n-m}-(a^{n-m})^{2}=a^{n-m-1}(a-a^{n-m+1}) \in \mathcal{A}^{qnil}$.
Thus, in view of Lemma \ref{lem3} we obtain $a^{n-m}$ is gs-Drazin invertible.
Let $x$ be the gs-Drazin inverse of $a^{n-m}$, i.e.
$$xa^{n-m}x=x,\ xa^{n-m}=a^{n-m}x \ \mbox{and}\ a^{n-m}-a^{n-m}x\in \mathcal{A}^{qnil}.$$

Define $y=a^{n-m-1}x$.
Next, we prove that $a\in \mathcal{A}^{g\pi H}$ and $a^{g\pi H}=y$.
Observe the fact that $x\in \mbox{comm}^{2}(a^{n-m})$. So, $xa=ax$.
Then, we have $ay=ya$ and $yay=y$.
It is clear that
$$a-a^{(n-m)+2}y=(a-a^{2n-2m+1})+(a^{2n-2m+1}-a^{2n-2m+1}x).$$
Since $$a-a^{2n-2m+1}=(a-a^{n-m+1})(1+a^{n-m})\in \mathcal{A}^{qnil}$$
and
$$a^{2n-2m+1}-a^{2n-2m+1}x=a^{n-m+1}(a^{n-m}-a^{n-m}x)\in \mathcal{A}^{qnil},$$
from Lemma \ref{lem1}(2) it follows that $a-a^{(n-m)+2}y\in \mathcal{A}^{qnil}$.
This completes the proof.
\qed

\medskip
Applying Theorem \ref{theorem4}, we get the following result.
\begin{cor}\label{cor1}
Let $a\in \mathcal{A}$ and $k\in \mathbb{N}$.
Then, $a\in \mathcal{A}^{g\pi H}$ if and only if $a^{k}\in \mathcal{A}^{g\pi H}$.
\end{cor}
\proof
Suppose that $a\in \mathcal{A}^{g\pi H}$.
By Theorem \ref{theorem4}(1)(2), we get $a-a^{n+1}\in \mathcal{A}^{qnil}$, for some $n\in \mathbb{N}$.
Then, we deduce that
$$a^{k}-(a^{k})^{n+1}=a^{k}-(a^{n+1})^{k}=\left(a-a^{n+1}\right)\sum\limits_{i=0}^{k-1}a^{ni+k-1}\in \mathcal{A}^{qnil},$$
which implies $a^{k}\in \mathcal{A}^{g\pi H}$.

On the contrary, we have $a^{k}-(a^{k})^{m+1}\in \mathcal{A}^{qnil}$ for some $m\in \mathbb{N}$.
So, $a^{k}-a^{km+k}\in \mathcal{A}^{qnil}$.
Evidently, $k\neq km+k$.
According to Theorem \ref{theorem4}(1)(3), it follows that $a\in \mathcal{A}^{g\pi H}$.
\qed

\medskip
Now, we are in the position to give the equivalent definition for the g$\pi$-Hirano inverse in a Banach algebra.
\begin{thm}\label{thm2}
Let $a,x\in \mathcal{A}$.  Then the following are equivalent:

\emph {(1)} $a\in \mathcal{A}^{g\pi H}$ and $a^{g\pi H}=x$;

\emph {(2)} $xax=x$, $xa=ax$ and $a^{n}-ax\in \mathcal{A}^{qnil}$, for some $n\in \mathbb{N}$;

\emph {(3)} $xax=x$, $xa=ax$ and $a^{n}-a^{m}x\in \mathcal{A}^{qnil}$, for some $m,n\in \mathbb{N}$ such that $m-n\neq 1$.
\end{thm}
\proof
(1) $\Rightarrow$ (2).
It is clear that
 $$xax=x, xa=ax \ \mbox{and}\  a-a^{n+2}x\in \mathcal{A}^{qnil},$$
for some $n\in \mathbb{N}$.
According to the proof of the implication (1) $\Rightarrow$ (2) of Theorem \ref{theorem4} and Proposition \ref{prop1}, we  see that
$a-a^{n+1}\in \mathcal{A}^{qnil}$ and $a-a^{2}x\in \mathcal{A}^{qnil}$.
Thus, by Lemma \ref{lem1} we obtain
$$a^{n}-ax=a^{n-1}(a-a^{2}x)-(a-a^{n+1})x\in \mathcal{A}^{qnil}.$$

(2) $\Rightarrow$ (3). It is trivial.

(3) $\Rightarrow$ (1). Using item (3),
we get
$$(a-a^{2}x)^{n}=a^{n}(1-ax)=(a^{n}-a^{m}x)(1-ax) \in \mathcal{A}^{qnil},$$
i.e. $a-a^{2}x\in \mathcal{A}^{qnil}$.
Since $m-n\neq 1$, then we can consider  the following two cases.

Case 1: Assume that $m-n\geq 2$.
Then,
$$a^{n}-a^{m-1}=(a^{n}-a^{m}x)-a^{m-2}(a-a^{2}x) \in \mathcal{A}^{qnil}.$$
Thus, we deduce that
$$(a-a^{m-n})^{n}=(a^{n}-a^{m-1})(1-a^{m-n-1})^{n-1} \in \mathcal{A}^{qnil},$$
i.e. $a-a^{m-n}\in \mathcal{A}^{nil}$.
So,
$$a-a^{(m-n-1)+2}x=(a-a^{2}x)+(a-a^{m-n})ax \in \mathcal{A}^{qnil}.$$

Case 2: Assume that $n-m\geq 0$.
By the hypotheses $a^{n}-a^{m}x\in \mathcal{A}^{qnil}$ and $ax=xa$,
we conclude that
$$(a^{n-m+1}-ax)^{m}=(a^{n}-a^{m}x)(a^{n-m}-x)^{m-1} \in \mathcal{A}^{qnil}.$$
Hence, we get $a^{n-m+1}-ax\in \mathcal{A}^{qnil}$, which yields
$$x-a^{n-m+1}x=-x(a^{n-m+1}-ax)\in \mathcal{A}^{qnil}.$$
Then,
$$a-a^{(n-m+1)+2}x=(a-a^{2}x)+(x-a^{n-m+1}x)a^{2}\in \mathcal{A}^{qnil}.$$

Therefore, by these two cases we obtain $a\in \mathcal{A}^{g\pi H}$ and $a^{g\pi H}=x$.
\qed

\medskip
From Theorem \ref{thm2} we can see the relationship between  the g$\pi$-Hirano inverse and the generalized $n$-strong Drazin inverse.
Also, the equivalent definitions of g-Drazin inverse are obtained as follows.
\begin{rem}
(i) In \cite{Dijana19}, Mosi\'{c} introduced the definition of the generalized $n$-strong Drazin inverse (or gns-Drazin inverse),
where $n$ is a fixed positive integer.
By Theorem \ref{thm2}(1)(2), we can see that if $a\in \mathcal{A}$ is  gns-Drazin invertible then $a$ is
g$\pi$-Hirano invertible.
Conversely, the g$\pi$-Hirano  inverse is a kind of the gns-Drazin inverse.

(ii) In item (2) of Theorem \ref{thm2}, for the case $m-n=1$, we have $a\in  \mathcal{A}^{d}$ with $a^{d}=x$ if and only if
\begin{center}
$xax=x, xa=ax \ \mbox{and}\  a^{n}-a^{n+1}x\in \mathcal{A}^{qnil},$
\end{center}
for some $n\in \mathbb{N}$.

(iii) By Proposition \ref{prop1} and Theorem \ref{thm2},
it is evident that
$a\in \mathcal{A}^{d}$ and $a^{d}=x$ if and only if
\begin{center}
 $ xax=x$, $ax=xa$ $\mathrm{and}$ $a-a^{n}x\in \mathcal{A}^{qnil}$,
 \end{center}
for some $n\in \mathbb{N}$,
if and only if
\begin{center}
$ xax=x$, $ax=xa$ $\mathrm{and}$ $a^{m}-a^{n}x\in \mathcal{A}^{qnil}$,
\end{center}
for some $m,n\in \mathbb{N}$.
\end{rem}

\section{Additive results on the g$\pi$-Hirano invertibility}\label{a}
Let $a,b\in \mathcal{A}$ and $k\in \mathbb{N}$.
Then, the elements $a,b$ are said to satisfy the ``$k\star$" condition if
$$ab\prod\limits_{i=1}^{k}\alpha_{i}=0,\ \ \mbox{for\ any}\ \alpha_{1}, \alpha_{2}, \cdots, \alpha_{k} \in \{a,b\}.$$
Obviously, if $a,b$ satisfy the ``$k\star$" condition, then $a,b$ satisfy the ``$(k+1)\star$" condition,
but $b,a$ do not satisfy the ``$k\star$" condition in general.
Note that if $ab=0$ then  $a,b$ satisfy the ``$k\star$" condition, for any $k\in \mathbb{N}$.
Also, for $k=1,2,3$, the ``$k\star$" condition become the following special cases, respectively.

(1) $aba=ab^{2}=0$;

(2) $abab=aba^{2}=ab^{2}a=ab^{3}=0$;

(3) $ababa=abab^{2}=aba^{3}=aba^{2}b=ab^{2}a^{2}=ab^{2}ab=ab^{3}a=ab^{4}=0$.

\medskip
The ``$k\star$" condition was introduced by Cvetkovi\'{c}-Ili\'{c} \cite{Cvetkovic-Ilic13}.
For two Drazin invertible elements $a,b$ in a ring,
the author studied the sufficient condition for the Drazin invertibility of the sum $a+b$ under  the ``$k\star$" condition.
Motivated by this, in this section we will consider the equivalence of the g$\pi$-Hirano invertibility between  the elements $a,b$ and the sum $a+b$ under the ``$k\star$" condition in a Banach algebra.

\medskip
We begin with the following crucial lemma.
\begin{lem}\label{lem5}
Let $a,b\in \mathcal{A}$  and $k\in \mathbb{N}$.
If $a,b$ satisfy the \emph{``}$k\star$\emph{"} condition, then
$$a, b\in \mathcal{A}^{qnil} \Longleftrightarrow a+b \in \mathcal{A}^{qnil}.$$
\end{lem}
\proof

Suppose that $a, b\in \mathcal{A}^{qnil}$.
Note that
$$\begin{array}{ccl}
(a+b)^{k+2}&=&a^{2}(a+b)^{k}+\left(ba(a+b)^{k}+b^{2}(a+b)^{k}\right)\\
&:=&x_{1}+\left(x_{2}+x_{3}\right).
\end{array}$$
Since $a,b\in \mathcal{A}$ satisfy the \emph{``}$k\star$\emph{"} condition,
we have $x_{1}(x_{2}+x_{3})=0$, $x_{2}x_{3}=0$ and $x_{2}^{2}=0$.
Obviously,
$x_{1}=a^{k+2}+\sum\limits_{i=1}^{2^{k}-1}p_{i}abq_{i}$ for suitable $p_{i}, q_{i}\in \mathcal{A}$, where $i\in \overline{1,2^{k-1}}$.
Observe that $\left(\sum\limits_{i=1}^{2^{k}-1}p_{i}abq_{i}\right)^{2}=0$, $a^{k+2}\in \mathcal{A}^{qnil}$ and
$\left(\sum\limits_{i=1}^{2^{k}-1}p_{i}abq_{i}\right)a^{k+2}=0$.
Thus, in view of Lemma \ref{lem121} we obtain $x_{1}\in \mathcal{A}^{qnil}$.
Similarly, we can get $x_{3}\in \mathcal{A}^{qnil}$.
So, we have $(a+b)^{k+2}\in \mathcal{A}^{qnil}$, i.e.
$a+b \in \mathcal{A}^{qnil}$.

Conversely, let us suppose that $a+b \in \mathcal{A}^{qnil}$.
Applying the \emph{``}$k\star$\emph{"}  condition, we have  the following equations:
$$a(a+b)^{m}a^{k}=a^{m+k+1} \ \mbox{and}\  b(a+b)^{m}b^{k+1}=b^{m+k+2},\ \mbox{for any} \ m\in \mathbb{N}.$$
Therefore, we get
$$\|a^{m+k+1}\|=\|a(a+b)^{m}a^{k}\|\leq \|a\|^{k+1}\|\|(a+b)^{m}\|,$$
which together with $a+b \in \mathcal{A}^{qnil}$  imply that
$$r(a)=\lim\limits_{m\to\infty}\left(\|a^{m+k+1}\|^{\frac{1}{m}}\right)^{\frac{m}{m+k+1}}=\lim\limits_{m\to\infty}\|a^{m+k+1}\|^{\frac{1}{m}}\\
\leq \lim\limits_{m\to\infty}\|a\|^{\frac{k+1}{m}}\|(a+b)^{m}\|^{\frac{1}{m}}=0.$$
Therefore, $a \in \mathcal{A}^{qnil}$.
Similarly, we can verify $b \in \mathcal{A}^{qnil}$.
\qed

\medskip
Now, we give the relationship between  the g$\pi$-Hirano invertibility of $a$, $b$ and that of the sum $a+b$ under the ``$k\star$\emph{"} condition in a Banach algebra as follows.
\begin{thm}\label{thm1}
Let $a,b\in \mathcal{A}$  and $k\in \mathbb{N}$.
If $a,b$ satisfy the \emph{``}$k\star$\emph{"} condition, then
$$a, b\in \mathcal{A}^{g\pi H} \Longleftrightarrow  a+b \in \mathcal{A}^{g\pi H}.$$
\end{thm}
\proof
Suppose that $a,b\in \mathcal{A}^{g\pi H}$.
So, there exist $m_{1}$, $m_{2}\in \mathbb{N}$ such that
$a-a^{m_{1}+1}\in \mathcal{A}^{qnil}$ and $b-b^{m_{2}+1}\in \mathcal{A}^{qnil}$.
Take $m=km_{1}m_{2}$.
Then, we get $m\geq k$ and $$a-a^{m+1}=\left(a-a^{m_{1}+1}\right)\left(1+a^{m_{1}}+a^{2m_{1}}+\cdots+a^{(km_{2}-1)m_{1}}\right)\in \mathcal{A}^{qnil}.$$
Similarly, $b-b^{m+1}\in \mathcal{A}^{qnil}$.
Note that
$$\begin{array}{ccl}
x&:=&\left(a+b\right)-\left(a+b\right)^{m+1}\\
&=&\left(a-a^{m+1}\right)+\left(b-b^{m+1}\right)+\left(\sum\limits_{i}s_{i}abt_{i}+\sum\limits_{i}u_{i}bav_{i}\right)\\
&:=&x_{1}+x_{2}+x_{3},
\end{array}$$
where $s_{i}, t_{i}, u_{i}, v_{i}\in \mathcal{A}$.
Since $a,b$ satisfy the \emph{``}$k\star$\emph{"} condition,
by computation we  conclude that $$\left(\sum\limits_{i}s_{i}abt_{i}\right)^{2}=0, \ \left(\sum\limits_{i}u_{i}bav_{i}\right)^{3}=0\
\mbox{and} \ \left(\sum\limits_{i}s_{i}abt_{i}\right)\left(\sum\limits_{i}u_{i}bav_{i}\right)=0.$$
Hence,  $x_{3}\in \mathcal{A}^{qnil}$.
Clearly, $x_{3}$ and $x_{2}$ satisfy the \emph{``}$k\star$\emph{"} condition.
Then, by Lemma \ref{lem5} it follows that  $x_{2}+x_{3}\in \mathcal{A}^{qnil}$.
In addition, note that
$x_{1}$ and $x_{2}+x_{3}$ also satisfy the \emph{``}$k\star$\emph{"} condition.
Applying Lemma \ref{lem5} again, we derive
$x\in \mathcal{A}^{qnil}$, which yields $a+b \in \mathcal{A}^{g\pi H}$.

Similar to the statements above,
 we can prove the sufficiency of this theorem in terms of Lemma \ref{lem5}.
\qed

\medskip
The following corollary can be directly obtained from Theorem \ref{thm1}.
\begin{cor}
Let $a,b\in \mathcal{A}$. If $ab=0$, then
$$a, b\in \mathcal{A}^{g\pi H}  \Longleftrightarrow a+b \in \mathcal{A}^{g\pi H}.$$
\end{cor}

Following the same strategy as in the proof of Theorem \ref{thm1}, we have
\begin{thm}
Let $a,b\in \mathcal{A}$ and $k\in \mathbb{N}$. If $\left(\prod\limits_{i=1}^{k}\alpha_{i}\right)ab=0$ for any $\alpha_{1}, \alpha_{2}, \cdots, \alpha_{k} \in \{a,b\}$,
then
$$a, b\in \mathcal{A}^{g\pi H} \Longleftrightarrow  a+b \in \mathcal{A}^{g\pi H}.$$
\end{thm}

\medskip

\begin{rem}
Let us compare the g-Drazin invertibility with g$\pi$-Hirano invertibility for the sum $a+b$ under the condition $ab=0$.
It is well known that the following holds: for $a,b \in \mathcal{A}$ satisfying $ab=0$, then
$$a, b\in \mathcal{A}^{d} \Longrightarrow a+b \in \mathcal{A}^{d}.$$
However, we do not know whether the converse of the above implication holds or not.
If we consider the g-Drazin invertibility under the hypothesis $ab=ba=0$,
then we have the following result.
\end{rem}

\begin{thm}
Let $a,b\in \mathcal{A}$. If $ab=ba=0$,
then
$$a, b\in \mathcal{A}^{d} \Longleftrightarrow  a+b \in \mathcal{A}^{d}.$$
\end{thm}
\proof
The necessity follows directly by \cite[Theorem 5.7]{{Koliha96}}.

Now, suppose that  $a+b \in \mathcal{A}^{d}$.
We only need to prove  $a\in \mathcal{A}^{d}$ by the symmetry of $a,b$.
Let $x_{1}=a^{2}(a+b)^{d}$ and $x_{2}=a-a^{2}(a+b)^{d}$.
Then,  $a=x_{1}+x_{2}$.
Since $ab=ba=0$, then we have
$a(a+b)=(a+b)a$, which yields $a(a+b)^{d}=(a+b)^{d}a$.
Combining the following equalities
$$a\left((a+b)^{d}\right)^{2}=a(a+b)\left((a+b)^{d}\right)^{3}=a^{2}\left((a+b)^{d}\right)^{3}=\cdots=a^{n}\left((a+b)^{d}\right)^{n+1}$$
for any $n\in \mathbb{N}$,
we can verify that $x_{1}\in \mathcal{A}^{d}$ and $x_{1}^{d}=a\left((a+b)^{d}\right)^{2}$.
Using   $$\left((a+b)-(a+b)^{2}(a+b)^{d}\right)x_{2}=x_{2}^{2}=x_{2}\left((a+b)-(a+b)^{2}(a+b)^{d}\right),$$
we get $x_{2}\in \mathcal{A}^{qnil}$ by Lemma \ref{lem1}(1) and Lemma \ref{lem18}.
Note that $x_{1}x_{2}=x_{1}x_{2}=0$, so we obtain  $a\in \mathcal{A}^{d}$ according to the necessity of this theorem. \qed

\medskip
In order to continue considering the topic on  the g$\pi$-Hirano invertibility of the sum $a+b$,
we need to prepare the following.

At the first we present a new condition.
Namely, for $a,b\in \mathcal{A}$ and $k\in \mathbb{N}$,
we say that $a,b$  satisfy the ``$k\ast$" condition, i.e.
$$\left(\prod\limits_{i=1}^{k}\alpha_{i}\right)ab=\left(\prod\limits_{i=1}^{k}\alpha_{i}\right)ba,\ \mbox{for\ any}\ \alpha_{1}, \alpha_{2}, \cdots, \alpha_{k} \in \{a,b\}.$$
We can see that if $a,b$ satisfy the ``$k\ast$" condition then $a,b$ satisfy the ``$(k+1)\ast$" condition.
In addition, the ``$k\ast$" condition contains the following specializations:

(1) $ab=ba$;

(2) $a^{2}b=aba$ and $b^{2}a=bab$;

(3) $a^{3}b=a^{2}ba$, $ba^{2}b=(ba)^{2}$, $ab^{2}a=(ab)^{2}$ and $b^{2}ab=b^{3}a$.

\medskip

Let $e \in \mathcal{A}$   be an idempotent ($e^{2}=e$). Then we can represent element  $a \in \mathcal{A}$ as

           $$a=\left(\begin{matrix}
                    a_{1}&a_{3} \\
                    a_{4}&a_{2} \\
             \end {matrix}\right) _{e},$$
where $a_{1}=eae$, $a_{2}=(1-e)a(1-e)$,  $a_{3}=ea(1-e)$ and  $a_{4}=(1-e)ae$.

\medskip

In what follows, by $ \mathcal{A}_{1}$, $\mathcal{A}_{2}$ we denote the algebra $p\mathcal{A }p$, $(1-p)\mathcal{A}(1-p) $, where $p^{2}=p \in \mathcal{A}$, respectively.
The following lemmas play an important role in the sequel.
\begin{lem}\label{lem7}
Let $x,y\in \mathcal{A}$ and $p^{2}=p\in \mathcal{A}$.
If $x$ and $y$ have the representations
$$
x= \left(\begin{matrix}
a&c \\
0&b \\
\end {matrix}\right) _{p}\ \ and\ \
y= \left(\begin{matrix}
b&0\\
c&a\\
\end {matrix}\right) _{1-p},
$$
then
\emph{(1)} $a\in \mathcal{A}_{1}^{qnil}\ \mbox{and}\ b\in \mathcal{A}_{2}^{qnil}  \Longleftrightarrow x\in \mathcal{A}^{qnil}\ (resp. \ y\in \mathcal{A}^{qnil})$;

\ \ \emph{(2)} $a\in \mathcal{A}_{1}^{g\pi H}\ and \ b\in \mathcal{A}_{2}^{g\pi H}  \Longleftrightarrow x\in \mathcal{A}^{g\pi H}\ (resp. \ y\in \mathcal{A}^{g\pi H}).$
\end{lem}
\proof
(1) Assume that $a\in \mathcal{A}_{1}^{qnil}\ \mbox{and}\ b\in \mathcal{A}_{2}^{qnil}$.
Since $\sigma_{\mathcal{A}}(x)\subseteq \sigma_{\mathcal{A}_{1}}(a)\cup \sigma_{\mathcal{A}_{2}}(b)$,
then we get $\sigma_{\mathcal{A}}(x)=\{0\}$, i.e. $x\in \mathcal{A}^{qnil}$.

On the converse, note that $(1-p)xp=0$, i.e. $pxp=xp$.
Then, by induction we obtain that $a^{m}=(pxp)^{m}=px^{m}p$ for any $m\in \mathbb{N}$.
Using the condition $x\in \mathcal{A}^{qnil}$,
we conclude  $$r(a)=\lim\limits_{m\to\infty}\|a^{m}\|^{\frac{1}{m}}=\lim\limits_{m\to\infty}\|px^{m}p\|^{\frac{1}{m}}
\leq \lim\limits_{m\to\infty}\|p\|^{\frac{1}{m}}\|x^{m}\|^{\frac{1}{m}}\|p\|^{\frac{1}{m}}=0.$$
So, $a\in \mathcal{A}_{1}^{qnil}$.
Also, by $\sigma_{\mathcal{A}_{2}}(b)\subseteq \sigma_{\mathcal{A}_{1}}(a)\cup \sigma_{\mathcal{A}}(x)$ it follows that $\sigma_{\mathcal{A}_{2}}(b)=\{0\}$,
i.e. $b\in \mathcal{A}_{2}^{qnil}$.

(2) For any $n\in \mathbb{N}$ we have
$$x-x^{n+1}= \left(\begin{matrix}
a&c\\
0&b\\
\end {matrix}\right)_{p}-
\left(\begin{matrix}
a&c\\
0&b\\
\end {matrix}\right)_{p}^{n+1}
= \left(\begin{matrix}
a-a^{n+1}&\triangle\\
0&b-b^{n+1}\\
\end {matrix}\right)_{p}.
$$
Suppose that $a\in \mathcal{A}_{1}^{g\pi H}$ and $b\in \mathcal{A}_{2}^{g\pi H}$.
Then, there exists $m\in \mathbb{N}$ such that
$a-a^{m+1}\in \mathcal{A}_{1}^{qnil}$
and $b-b^{m+1}\in \mathcal{A}_{2}^{qnil}$.
So, $x-x^{m+1}\in \mathcal{A}^{qnil}$ by item (1).
Therefore, we get $x\in \mathcal{A}^{g\pi H}$.
The sufficiency can be proved similarly.
\qed

\begin{rem}
Item (2) of Lemma \ref{lem7} is somewhat different from the g-Drazin inverse case, namely,
if $a\in \mathcal{A}_{1}^{d}$, then $ b\in \mathcal{A}_{2}^{d}$  if and only if $x\in \mathcal{A}^{d}$ (\cite[Theorem 2.3]{CK04}).
\end{rem}

\begin{lem}\label{lemm107}
Let $a,b\in \mathcal{A}$  and $k\in \mathbb{N}$.
If $a,b$ satisfy the \emph{``}$k\ast$\emph{"} condition, then

\emph{(1)} if $a\in \mathcal{A}^{qnil}$ \emph{(}or $b\in \mathcal{A}^{qnil}$\emph{)}, then $ab\in \mathcal{A}^{qnil}$;

\emph{(2)}  if $a\in \mathcal{A}^{qnil}$, then $b\in \mathcal{A}^{qnil}$ if and only if $a+b\in \mathcal{A}^{qnil}$;

\emph{(3)} if $a,b\in \mathcal{A}^{g\pi H}$, then $ab\in \mathcal{A}^{g\pi H}$;

\emph{(4)} if $a\in \mathcal{A}^{qnil}$,  $b\in \mathcal{A}^{g\pi H}$, then $a+b\in  \mathcal{A}^{g\pi H}$.
\end{lem}
\proof
(1) Since $a,b$ satisfy the \emph{``}$k\ast$\emph{"} condition, then
we conclude that $(ab)^{n+k}=(ab)^{k}a^{n}b^{n}$ for any $n\in \mathbb{N}$.
Applying the hypothesis $a\in \mathcal{A}^{qnil}$ or $b\in \mathcal{A}^{qnil}$, we obtain
$$r(ab)=\lim\limits_{n\to\infty}\|(ab)^{n+k}\|^{\frac{1}{n}}\leq \lim\limits_{n\to\infty}\|(ab)^{k}\|^{\frac{1}{n}} \lim\limits_{n\to\infty}\|a^{n}\|^{\frac{1}{n}}\lim\limits_{n\to\infty}\|b^{n}\|^{\frac{1}{n}}=0,$$
which implies $ab\in \mathcal{A}^{qnil}$.

(2) Suppose that $a, b\in \mathcal{A}^{qnil}$. Note that $(a+b)^{k+1}=(a+b)^{k}a+(a+b)^{k}b$.
Let $c=(a+b)^{k}a$ and $d=(a+b)^{k}b$.
From the \emph{``}$k\ast$\emph{"} condition, we get $cd=dc$ and
$c^{n}=(a+b)^{kn}a^{n}$ for any $n\in \mathbb{N}$.
So,  $$r(c)=\lim\limits_{n\to\infty}\|c^{n}\|^{\frac{1}{n}}\leq \|a+b\|^{k}\lim\limits_{n\to\infty}\|a^{n}\|^{\frac{1}{n}}=0,$$
which means $c\in \mathcal{A}^{qnil}$.
Similarly, we have $d\in \mathcal{A}^{qnil}$.
Applying Lemma \ref{lem1}(2), it follows that $(a+b)^{k+1}\in \mathcal{A}^{qnil}$, i.e. $a+b\in \mathcal{A}^{qnil}$.

To prove the converse, let us suppose that $a, a+b\in \mathcal{A}^{qnil}$.
Obviously, $-a,a+b$ satisfy the \emph{``}$k\ast$\emph{"} condition.
Then,  we deduce that $b=-a+(a+b)\in \mathcal{A}^{qnil}$ by the proof of the necessity of item (2).

(3) In view of the condition $a,b\in \mathcal{A}^{g\pi H}$, we obtain
$a-a^{m+1}\in \mathcal{A}^{qnil}$ and  $b-b^{m+1}\in \mathcal{A}^{qnil}$ for some $m\in \mathbb{N}$.
By the \emph{``}$k\ast$\emph{"} condition, we get
$$(ab)^{k+1}-(ab)^{m+k+1}=(ab)^{k}(a-a^{m+1})b+(ab)^{k}a^{m+1}(b-b^{m+1}).$$
Setting $s=(ab)^{k}(a-a^{m+1})b$ and  $t=(ab)^{k}a^{m+1}(b-b^{m+1})$.
Applying the \emph{``}$k\ast$\emph{"} condition again, we obtain
$s^{n}=(ab)^{kn}(a-a^{m+1})^{n}b^{n}$ for any $n\in \mathbb{N}$, which implies
$$r(s)=\lim\limits_{n\to\infty}\|s^{n}\|^{\frac{1}{n}}\leq\|ab\|^{k}\|b\|\lim\limits_{n\to\infty}\|(a-a^{m+1})^{n}\|^{\frac{1}{n}}=0.$$
So, $s\in \mathcal{A}^{qnil}$.
Similarly, $t\in \mathcal{A}^{qnil}$.
Note that $st=ts$.
Therefore, $s+t\in \mathcal{A}^{qnil}$,
which gives $ab\in \mathcal{A}^{g\pi H}$ by Theorem \ref{theorem4}(3).

(4) Let $p=bb^{g\pi H}$.
Now, we  consider the matrix representations of $a$ and $b$ relative to the idempotent $p$:
$$a= \left(\begin{matrix}
a_{1}&a_{3}\\
a_{4}&a_{2}\\
\end {matrix}\right)_{p}  \ \mbox{and}\ \
 b= \left(\begin{matrix}
 b_{1}&0\\
 0&b_{2}\\
 \end {matrix}\right) _{p}.$$
Obviously, $b_{1}\in  \mathcal{A}_{1}^{-1}$ with $(b_{1})_{\mathcal{A}_{1}}^{-1}=b^{g\pi H}$.
Also, by item (3) we get $b_{1}=b(bb^{g\pi H})\in \mathcal{A}_{1}^{g\pi H}$.
From Proposition \ref{prop1}, it follows that  $b_{2}=b-b^{2}b^{g\pi H}=b-b^{2}b^{d}\in \mathcal{A}_{2}^{qnil}$.

Note that
$$ \left(\begin{matrix}
b_{1}^{k}a_{1}b_{1}&b_{1}^{k}a_{3}b_{2}\\
b_{2}^{k}a_{4}b_{1}&b_{2}^{k}a_{2}b_{2}\\
\end {matrix}\right)_{p}=b^{k}ab=b^{k+1}a=\left(\begin{matrix}
b_{1}^{k+1}a_{1}&b_{1}^{k+1}a_{3}\\
b_{2}^{k+1}a_{4}&b_{2}^{k+1}a_{2}\\
\end {matrix}\right)_{p}.$$
Thus, we have  $b_{1}^{k}a_{3}b_{2}=b_{1}^{k+1}a_{3}$, so, $a_{3}=b_{1}^{-1}a_{3}b_{2}$,
which implies $a_{3}=b_{1}^{-n}a_{3}b_{2}^{n}$ for any $n\in \mathbb{N}$.
Since $b_{2}\in \mathcal{A}_{2}^{qnil}$,
we have $$\lim\limits_{n\to\infty}\|a_{3}\|^{\frac{1}{n}}\leqslant \|b_{1}^{-1}\|\lim\limits_{n\to\infty}\|a_{3}\|^{\frac{1}{n}}\lim\limits_{n\to\infty}\|b_{2}^{n}\|^{\frac{1}{n}}=0.$$
Hence, $a_{3}=0$.
In addition, it is easy to see that  $a_{1}b_{1}=b_{1}a_{1}$ and $a_{2}, b_{2}$ satisfy the \emph{``}$k\ast$\emph{"} condition.

Now,  we have
$$a= \left(\begin{matrix}
 a_{1}&0\\
 a_{4}&a_{2}\\
 \end {matrix}\right)_{p}\ \mbox{and}\ \
 a+b=\left(\begin{matrix}
           a_{1}+b_{1}&0\\
           a_{4}&a_{2}+b_{2}\\
 \end {matrix}\right)_{p}.$$
From the condition $a\in \mathcal{A}^{qnil}$, it follows that $a_{1}\in \mathcal{A}_{1}^{qnil}$ and  $a_{2}\in \mathcal{A}_{2}^{qnil}$ by Lemma \ref{lem7}(1).
Applying item (2), we can obtain $a_{2}+b_{2}\in \mathcal{A}_{2}^{qnil}$,
which implies  $a_{2}+b_{2}\in \mathcal{A}_{2}^{g\pi H}$.
Note that $(p+b_{1}^{-1}a_{1})-(p+b_{1}^{-1}a_{1})^{2}=-a_{1}(b_{1}^{-1}+a_{1}b_{1}^{-2})\in \mathcal{A}_{1}^{qnil}$ by Lemma \ref{lem1}(1).
Hence, we conclude $p+b_{1}^{-1}a_{1}\in \mathcal{A}_{1}^{g\pi H}$.
Then, in view of item (3), we obtain $a_{1}+b_{1}=b_{1}(p+b_{1}^{-1}a_{1})\in \mathcal{A}_{1}^{g\pi H}$.
Finally, by Lemma \ref{lem7}(2) we deduce $a+b\in  \mathcal{A}^{g\pi H}$.
\qed

\medskip
Now, we present the equivalent characterization for the g$\pi$-Hirano invertibility of the sum $a+b$ under the \emph{``}$k\ast$\emph{"} condition.
\begin{thm}\label{thm88}
Let $a,b\in \mathcal{A}^{g\pi H}$  and $k\in \mathbb{N}$.
If $a,b$ satisfy the \emph{``}$k\ast$\emph{"} condition, then
$$1+a^{g\pi H}b\in \mathcal{A}^{g\pi H} \Longleftrightarrow a+b\in  \mathcal{A}^{g\pi H}.$$
\end{thm}
\proof
Let $p=aa^{g\pi H}$. Then,
as in the proof of Lemma \ref{lemm107}(4) we have
$$a= \left(\begin{matrix}
a_{1}&0\\
0&a_{2}\\
\end {matrix}\right)_{p}  \ \mbox{and}\ \
 b= \left(\begin{matrix}
 b_{1}&0\\
 b_{4}&b_{2}\\
 \end {matrix}\right)_{p},$$
 where $a_{1}\in  \mathcal{A}_{1}^{-1}\cap \mathcal{A}_{1}^{g\pi H}$ and $a_{2}\in \mathcal{A}_{2}^{qnil}$.
In addition, we have  $a_{1}b_{1}=b_{1}a_{1}$, and $a_{2}, b_{2}$ satisfy the \emph{``}$k\ast$\emph{"} condition.
Using  $b\in \mathcal{A}^{g\pi H}$, we get $b_{1}\in \mathcal{A}_{1}^{g\pi H}$
and $b_{2}\in \mathcal{A}_{2}^{g\pi H}$ by Lemma \ref{lem7}(2).

Note that
$$1+a^{g\pi H}b=\left(\begin{matrix}
p+a_{1}^{-1}b_{1}&0\\
0&1-p\\
\end {matrix}\right)_{p}\ \mbox{and}\ \
 a+b=\left(\begin{matrix}
           a_{1}+b_{1}&0\\
           b_{4}&a_{2}+b_{2}\\
 \end {matrix}\right)_{p}.$$
By Lemma \ref{lemm107}(4), we have $a_{2}+b_{2}\in \mathcal{A}_{2}^{g\pi H}$.
From Lemma \ref{lem7}(2), we claim that  $a+b\in  \mathcal{A}^{g\pi H}$ if and only if
$a_{1}+b_{1}\in  \mathcal{A}_{1}^{g\pi H}$,
and $1+a^{g\pi H}b\in \mathcal{A}^{g\pi H}$ if and only if
$p+a_{1}^{-1}b_{1}\in  \mathcal{A}_{1}^{g\pi H}$.

Next, we only need to show that $p+a_{1}^{-1}b_{1}\in  \mathcal{A}_{1}^{g\pi H}$ is equivalent to $a_{1}+b_{1}\in  \mathcal{A}_{1}^{g\pi H}$.
If $p+a_{1}^{-1}b_{1}\in  \mathcal{A}_{1}^{g\pi H}$,
then $a_{1}+b_{1}=a_{1}(p+a_{1}^{-1}b_{1})\in  \mathcal{A}_{1}^{g\pi H}$ by Lemma \ref{lemm107}(3).
On the contrary, let us suppose that $a_{1}+b_{1}\in  \mathcal{A}_{1}^{g\pi H}$.
In view of the hypothesis $a\in \mathcal{A}^{g\pi H}$, we have $a-a^{m+1}\in \mathcal{A}^{qnil}$ for some $m\in \mathbb{N}$.
Then we get
$$a_{1}^{-1}-a_{1}^{-m-1}=a^{g\pi H}-(a^{g\pi H})^{m+1}=-(a^{g\pi H})^{m+2}(a-a^{m+1})\in \mathcal{A}^{qnil},$$
which implies $a_{1}^{-1}\in  \mathcal{A}_{1}^{g\pi H}$.
So, we conclude that $p+a_{1}^{-1}b_{1}=a_{1}^{-1}(a_{1}+b_{1}) \in \mathcal{A}_{1}^{g\pi H}$.
This completes the proof.
 \qed

\medskip
Dual to Theorem \ref{thm88}, we have the following the result.
\begin{thm}\label{thm898}
Let $a,b\in \mathcal{A}^{g\pi H}$ and $k\in \mathbb{N}$.
If $ab\prod\limits_{i=1}^{k}\alpha_{i}=ba\prod\limits_{i=1}^{k}\alpha_{i}$ for any $\alpha_{1}, \alpha_{2}, \cdots, \alpha_{k} \in \{a,b\}$, then
$$1+a^{g\pi H}b\in \mathcal{A}^{g\pi H} \Longleftrightarrow a+b\in  \mathcal{A}^{g\pi H}.$$
\end{thm}

Let us notice that the implications of the sufficiency  in Theorem \ref{thm88} and Theorem~
\ref{thm898} still hold without the assumption  $b\in \mathcal{A}^{g\pi H}$,
but if we remove this assumption then the  implications of the necessity are not valid in general.
This can be illustrated by the following example.
Let $\mathcal{A}=\mathbb{C}$, $a=0$ and $b=2$.
Then,  it is obvious that  $1+a^{g\pi H}b=1\in \mathcal{A}^{g\pi H}$.
However, $a+b=2 \notin  \mathcal{A}^{g\pi H}$.

\section{Anti-triangular matrices involving the g$\pi$-Hirano inverse}\label{a}
In this section, we mainly consider some sufficient and necessary conditions for anti-triangular matrices
$\left(\begin{matrix}
a&b\\
c&0\\
\end{matrix}\right)$ over Banach algebras to be
 g$\pi$-Hirano invertible.

The authors \cite{ZCM17} found the anti-triangular matrix
 $\left(\begin{matrix}
 1&1\\
 c&0\\
\end {matrix}\right)\in M_{2}(\mathcal{A})$ is Drazin invertible if and only if $c\in \mathcal{A}$ is Drazin invertible.
But, for the  g$\pi$-Hirano case, we do not have the corresponding result, which can be seen from the following example.
Let $\mathcal{A}=\mathbb{C}$ and $c=1$.
Obviously, $1\in \mathcal{A}^{g\pi H}$.
But, $\left(\begin{matrix}
 1&1\\
 1&0\\
\end {matrix}\right)\notin M_{2}(\mathcal{A})^{g\pi H}$,
since $\left(\begin{matrix}
 1&1\\
 1&0\\
\end {matrix}\right)-\left(\begin{matrix}
 1&1\\
 1&0\\
\end {matrix}\right)^{n+1}\notin M_{2}(\mathcal{A})^{qnil}$, for any $n\in \mathbb{N}$.
So, what is the equivalent conditions for $\left(\begin{matrix}
 1&1\\
 c&0\\
\end {matrix}\right)\in M_{2}(\mathcal{A})$
to be  g$\pi$-Hirano invertible?

At the beginning, we consider the equivalent conditions  for the matrix $\left(\begin{matrix}
 1&1\\
 c&0\\
\end {matrix}\right) \in M_{2}(\mathcal{A})$ to be g-Hirano invertible.
The definition of the g-Hirano inverse was introduced by Chen and Sheibani \cite{CS2019F}, namely
an element $a \in \mathcal{A}$ has g-Hirano inverse if there exists $x \in \mathcal{A}$ such that
$$ xax=x,\ \ ax=xa\ \   \mathrm{and} \ \ a^{2}-ax\in \mathcal{A}^{qnil}.$$
Clearly, by Theorem \ref{thm2} we see that the g-Hirano inverse is a subclass of the g$\pi$-Hirano inverse.
Denote by $\mathcal{A}^{gH}$ the set of all g-Hirano invertible elements in $\mathcal{A}$.
\begin{thm}\label{thm2021117}
Let $M= \left(\begin{matrix}
 1&1\\
 c&0\\
\end {matrix}\right)\in M_{2}(\mathcal{A})$. Then,
$$c\in \mathcal{A}^{qnil} \Longleftrightarrow M\in  M_{2}(\mathcal{A})^{gH}.$$
\end{thm}
\proof
From \cite[Theorem 2.4]{CS2019F} it follows that $M\in  M_{2}(\mathcal{A})^{gH}$ if and only if
$N:=M-M^{3}=-\left(\begin{matrix}
 2c&c\\
 c^{2}&c\\
\end {matrix}\right)\in M_{2}(\mathcal{A})^{qnil}$.
Therefore, we only need to prove that
$c\in \mathcal{A}^{qnil}$ is equivalent to $N \in M_{2}(\mathcal{A})^{qnil}$.

Suppose that $c\in \mathcal{A}^{qnil}$.
Then, $N=-\left(\begin{matrix}
 c&0\\
 0&c\\
\end {matrix}\right)\left(\begin{matrix}
 2&1\\
 c&1\\
\end {matrix}\right)\in M_{2}(\mathcal{A})^{qnil}$ by Lemma \ref{lem1}(1).

On the contrary, by $N \in M_{2}(\mathcal{A})^{qnil}$
we get that $\left(\begin{matrix}
 2c-\lambda&c\\
 c^{2}&c-\lambda\\
\end {matrix}\right)$
is invertible, for any $\lambda\in\mathbb{C}\backslash \{0\}$.
Since $$\left(\begin{matrix}
 1&-1\\
 0&1\\
\end {matrix}\right)
\left(\begin{matrix}
2c-\lambda&c\\
 c^{2}&c-\lambda\\
\end {matrix}\right)
=\left(\begin{matrix}
 -c^{2}+2c-\lambda&\lambda\\
 c^{2}&c-\lambda\\
\end {matrix}\right),$$
we deduce that $\left(\begin{matrix}
 -c^{2}+2c-\lambda&\lambda\\
 c^{2}&c-\lambda\\
\end {matrix}\right)$ is invertible for any $\lambda\in\mathbb{C}\backslash \{0\}$.
Hence, there exists $\left(\begin{matrix}
x&y\\
z&w\\
\end {matrix}\right)\in M_{2}(\mathcal{A})$
such that
$$\left(\begin{matrix}
 -c^{2}+2c-\lambda&\lambda\\
 c^{2}&c-\lambda\\
\end {matrix}\right)\left(\begin{matrix}
x&y\\
z&w\\
\end {matrix}\right)=\left(\begin{matrix}
1&0\\
0&1\\
\end {matrix}\right)$$
and
$$\left(\begin{matrix}
x&y\\
z&w\\
\end {matrix}\right)\left(\begin{matrix}
 -c^{2}+2c-\lambda&\lambda\\
 c^{2}&c-\lambda\\
\end {matrix}\right)=\left(\begin{matrix}
1&0\\
0&1\\
\end {matrix}\right).$$
So, we can obtain the following equations
\begin{equation}\label{eq1}
 (-c^{2}+2c-\lambda)y+\lambda w=0,
\end{equation}
\begin{equation}\label{eq2}
  c^{2}y+(c-\lambda)w=1,
\end{equation}
\begin{equation}\label{eq3}
  x(-c^{2}+2c-\lambda)+yc^{2}=1,
\end{equation}
\begin{equation}\label{eq4}
 \lambda x+y(a-\lambda)=0.
\end{equation}
By the equation (1), we have $w=\frac{1}{\lambda}\left(c^{2}-2c+\lambda\right)y$, which
together with (2) imply
$$1=c^{2}y+\frac{1}{\lambda}\left(c-\lambda\right)\left(c^{2}-2c+\lambda\right)y=\frac{1}{\lambda}\left(c^{3}-2c^{2}+3\lambda c-\lambda^{2}\right)y.$$
So, $y$ is left invertible.
Similarly, using (3) and (4) we conclude that $y$ is right invertible.
Therefore,
$y$ is invertible and $y^{-1}=\frac{1}{\lambda}\left(c^{3}-2c^{2}+3\lambda c-\lambda^{2}\right)$.
So, $c^{3}-2c^{2}+3\lambda c-\lambda^{2}$ is invertible.
Hence, $0\notin \sigma(c^{3}-2c^{2}+3\lambda c-\lambda^{2})$ for any $\lambda\in\mathbb{C}\backslash \{0\}$.

Now, assume that there exists $t\in\mathbb{C}\backslash \{0\}$ such that $t\in \sigma(c)$.
Then, we can find $\lambda_{0}\in\mathbb{C}\backslash \{0\}$ satisfying $t^{3}-2t^{2}+3\lambda_{0}t-\lambda_{0}^{2}=0$.
So, $0\in \sigma(c^{3}-2c^{2}+3\lambda_{0} c-\lambda_{0}^{2})$, which  contradicts with  $0\notin \sigma(c^{3}-2c^{2}+3\lambda c-\lambda^{2})$ for any $\lambda\in\mathbb{C}\backslash \{0\}$.
Hence, $\sigma(c)=\{0\}$, i.e. $c\in \mathcal{A}^{qnil}$.
\qed

\begin{rem}
By Theorem \ref{thm2021117}, we get $$c\in \mathcal{A}^{qnil} \Longrightarrow M= \left(\begin{matrix}
 1&1\\
 c&0\\
\end {matrix}\right)\in  M_{2}(\mathcal{A})^{g\pi H}.$$
However, in general  the converse of the above implication does not hold,  which can be seen from the following example:
\end{rem}
\begin{eg}
Let $\mathcal{A}=\mathbb{C}$ and $c=-1$.
Observe that
$\left(\begin{matrix}
1&1\\
-1&0\\
\end{matrix}\right)^{6}=\left(\begin{matrix}
 1&0\\
 0&1\\
\end{matrix}\right)$.
Therefore, we get that
$\left(\begin{matrix}
 1&1\\
 -1&0\\
\end {matrix}\right)\in  M_{2}(\mathcal{A})^{g\pi H}$ and
$\left(\begin{matrix}
 1&1\\
 -1&0\\
\end {matrix}\right)^{g\pi H}=\left(\begin{matrix}
 0&-1\\
 1&1\\
\end {matrix}\right)$.
However, $-1\notin \mathcal{A}^{qnil}$.
\end{eg}

Next, we will consider the g$\pi$-Hirano invertibility for the  anti-triangular matrix $\left(\begin{matrix}
 a&b\\
 c&0\\
\end {matrix}\right)$ over Banach algebras.   For future reference we state two lemmas as follows.
\begin{lem}\label{lem8}
Let $a,b \in \mathcal{A}$.
Then, $ab\in \mathcal{A}^{g\pi H}$ if and only if $ba\in \mathcal{A}^{g\pi H}$.
\end{lem}
\proof If $ab\in \mathcal{A}^{g\pi H}$, then $ab-(ab)^{m+1}\in \mathcal{A}^{qnil}$, for some $m\in\mathbb{N}$.
By induction, we have
$\left((ba)^{2}-(ba)^{m+2}\right)^{n}=b\left(ab-(ab)^{m+1}\right)^{n}(ab)^{n-1}a$, for any  $n\in \mathbb{N}$.
Thus,
$$\lim\limits_{n\to\infty}\|\left((ba)^{2}-(ba)^{m+2}\right)^{n}\|^{\frac{1}{n}}\leq\|a\|\|b\|\lim\limits_{n\to\infty}\|(ab-(ab)^{m+1})^{n}\|^{\frac{1}{n}}=0.$$
So, $(ba)^{2}-(ba)^{m+2}\in \mathcal{A}^{qnil}$, which means $ba\in \mathcal{A}^{g\pi H}$.
\qed

\begin{lem}\label{lem9}
Let $M= \left(\begin{matrix}
 a&c\\
 0&b\\
\end {matrix}\right) \left( or \left(\begin{matrix}
 a&0\\
 d&b\\
\end {matrix}\right)\right)\in M_{2}(\mathcal{A})$.
Then,

\emph{(1)} $a\in \mathcal{A}^{qnil}\ \mbox{and}\ b\in \mathcal{A}^{qnil} \Longleftrightarrow M\in M_{2}( \mathcal{A})^{qnil}$.

\emph{(2)} $a\in \mathcal{A}^{g\pi H}\ \mbox{and}\ b\in \mathcal{A}^{g\pi H} \Longleftrightarrow M\in M_{2}( \mathcal{A})^{g\pi H}$.
\end{lem}
\proof (1) Suppose that $M=\left(\begin{matrix}
 a&c\\
 0&b\\
\end {matrix}\right)\in M_{2}( \mathcal{A})^{qnil}$.
Let $P=\left(\begin{matrix}
  1&0\\
  0&0\\
  \end {matrix}\right)\in M_{2}(\mathcal{A})$.
Then we have the  following  matrix representation of $M$ relative to the idempotent $P$:
$$ M=
\left(\begin{matrix}
A&C\\
0& B\\
 \end {matrix}\right)_{P}, \ \mbox{where}\
 A=\left(\begin{matrix}
a&0\\
0&0\\
\end {matrix}\right),\
 B=\left(\begin{matrix}
0&0\\
0&b\\
\end {matrix}\right)\ \mbox{and}\
C=\left(\begin{matrix}
0&c\\
0&0\\
\end {matrix}\right).$$
By Lemma \ref{lem7}(1),
we obtain $A \in \left(PM_{2}(\mathcal{A})P\right)^{qnil}$ and $B\in \left((I-P)M_{2}(\mathcal{A})(I-P)\right)^{qnil}$,
i.e. $\sigma_{PM_{2}(\mathcal{A})P}(A)=\{0\}$ and $\sigma_{(I-P)M_{2}(\mathcal{A})(I-P)}(B)=\{0\}$.
Note that $\sigma_{PM_{2}(\mathcal{A})P}(A)\cup \{0\}=\sigma_{M_{2}(\mathcal{A})}(A)$.
So, $\sigma_{M_{2}(\mathcal{A})}(A)=\{0\}$, which implies that $\lambda I-A$ is invertible, for any $\lambda\neq 0$.
Hence, $\lambda 1-a$ is invertible, so $\sigma_{\mathcal{A}}(a)=\{0\}$.
Therefore, $a\in \mathcal{A}^{qnil}$.
Similarly, we can get $b\in \mathcal{A}^{qnil}$.
On the contrary, applying  $\sigma_{M_{2}(\mathcal{A})}(M)\subseteq \sigma_{\mathcal{A}}(a)\cup \sigma_{\mathcal{A}}(b)$
we deduce $M\in M_{2}( \mathcal{A})^{qnil}$.

\smallskip
(2) By item (1) and Theorem \ref{theorem4}, item (2) holds directly.
\qed

\medskip

Now, we are ready to present an existence criterion for the g$\pi$-Hirano inverse of  the anti-triangular matrix under the ``$k\star$" condition as follows.
\begin{thm}\label{thm117}
Let $M= \left(\begin{matrix}
 a&b\\
 c&0\\
\end {matrix}\right)\in M_{2}(\mathcal{A})$ and $k\in \mathbb{N}$.
If $a,bc$ satisfy the \emph{``}$k\star$\emph{"} condition, then
$$a, bc\in \mathcal{A}^{g\pi H} \Longleftrightarrow M\in  M_{2}(\mathcal{A})^{g\pi H}.$$
\end{thm}
\proof
Note that
$$M=
    \left(\begin{matrix}
    1&0\\
    0&c\\
    \end {matrix}\right)
    \left(\begin{matrix}
    a&b\\
    1&0\\
    \end {matrix}\right).$$
By Lemma \ref{lem8}, we have
   $$M\in M_{2}(\mathcal{A})^{g\pi H}  \Longleftrightarrow  N:=\left(\begin{matrix}
     a&b\\
    1&0\\
    \end {matrix}\right)
  \left(\begin{matrix}
    1&0\\
    0&c\\
    \end {matrix}\right)=
   \left(\begin{matrix}
    a&bc\\
    1&0\\
 \end {matrix}\right) \in M_{2}(\mathcal{A})^{g\pi H}.$$
Consider the following decomposition:
$$N^{2}=
    \left(\begin{matrix}
    a^{2}&0\\
    a&0\\
    \end {matrix}\right)+
    \left(\begin{matrix}
    bc&abc\\
    0&bc\\
    \end {matrix}\right):=N_{1}+N_{2}.$$
Since $a,bc$ satisfy the   \emph{``}$k\star$\emph{"} condition, so do $N_{1}$ and $N_{2}$.
Therefore, by using Corollary \ref {cor1}, Theorem \ref{thm1} and Lemma \ref{lem9}(2) we deduce that
$$N \in M_{2}(\mathcal{A})^{g\pi H}\Longleftrightarrow N^{2} \in M_{2}(\mathcal{A})^{g\pi H}
\Longleftrightarrow N_{1}, N_{2} \in M_{2}(\mathcal{A})^{g\pi H}\Longleftrightarrow a, bc\in \mathcal{A}^{g\pi H},$$
as required.
\qed

\medskip
It is easy to see that if the hypothesis $a,bc$ satisfy the \emph{``}$k\star$\emph{"} condition in Theorem \ref{thm117}  is replaced by
the hypothesis $bc, a$ satisfy the \emph{``}$k\star$\emph{"} condition then this theorem still holds.
So, we immediately obtain the  following corollary.

\begin{cor}\label{cor8}
Let $M= \left(\begin{matrix}
 a&b\\
 c&0\\
\end {matrix}\right)\in M_{2}(\mathcal{A})$.
If $abc=0$ $($or $bca=0$$)$, then
$$a, bc\in \mathcal{A}^{g\pi H} \Longleftrightarrow M\in  M_{2}(\mathcal{A})^{g\pi H}.$$
\end{cor}

Let us remark that in Corollary \ref{cor8} the condition $abc=0$ or $bca=0$ in general can not be substituted by $acb=0$ or $cab=0$,
which can be seen from the following example.
\begin{eg}
Let $\mathcal{A}=M_{2}(\mathbb{C})$
and  $M= \left(\begin{matrix}
 a&b\\
 c&0\\
\end {matrix}\right)\in M_{2}(\mathcal{A})$.
(1) If we choose $a=\left(\begin{matrix}
 1&0\\
 0&0\\
\end {matrix}\right)$,
$b=\left(\begin{matrix}
 0&1\\
 0&0\\
\end {matrix}\right)$ and
$c=\left(\begin{matrix}
 0&0\\
 1&0\\
\end {matrix}\right)$,
then $a, bc\in \mathcal{A}^{g\pi H}$ and $acb=0$.
However, $M-M^{n+1}\notin \mathcal{A}^{qnil}$  for any $n\in \mathbb{N}$,
so $M\notin  M_{2}(\mathcal{A})^{g\pi H}$.
(2) If we take  $a=\left(\begin{matrix}
 0&1\\
 1&0\\
\end {matrix}\right)$ and
$b=c=\left(\begin{matrix}
 1&0\\
 0&0\\
\end {matrix}\right)$,
then $a, bc\in \mathcal{A}^{g\pi H}$ and $cab=0$.
But, $M\notin  M_{2}(\mathcal{A})^{g\pi H}$.
\end{eg}

\medskip

\vspace{0.2cm} \noindent {\large\bf Acknowledgements}.
This research was supported by  NSF of Jiangsu Province (No. BK20200944),
Natural Science Foundation of Jiangsu Higher Education Institutions of China (No. 20KJB110001),
China Postdoctoral Science Foundation (No. 2020M671281).

\end{document}